\documentclass[12pt]{amsart}
\usepackage{amssymb, amsmath, amsfonts, amsthm, esint,graphicx}
\usepackage[abbrev]{amsrefs}
\usepackage[margin=1 in]{geometry}
\usepackage[colorlinks=true, linkcolor=blue]{hyperref}

\def \dd {\partial}

\def \eps {\varepsilon}

\DeclareMathOperator{\Ric}{Ric}

\DeclareMathOperator{\Hess}{Hess}

\DeclareMathOperator{\vol}{vol}

\DeclareMathOperator{\diver}{div}

\DeclareMathOperator{\SecFun}{II}

\title{First eigenvalue of the $p$-Laplacian on K\"ahler manifolds}
\author{Casey Blacker}
\address{Department of Mathematics\\
         University of California\\
         Santa Barbara, CA 93106}
\email{\href{mailto:clacker@ucsb.edu}{cblacker@ucsb.edu}}

\author{Shoo Seto}
\address{Department of Mathematics\\
         University of California\\
         Santa Barbara, CA 93106}
\email{\href{mailto:shoseto@ucsb.edu}{shoseto@ucsb.edu}}
\thanks{Partially support by Simons Travel Grant}
\date{}

\theoremstyle{definition} 
\newtheorem{definition}{Definition}[section]

\newtheorem{lemma}{Lemma}[section]
\newtheorem{remark}{Remark}[section]

\newtheorem{theorem}{Theorem}[section]

\newtheorem*{acknowledgements}{Acknowledgements}

\begin{document}
\maketitle
\begin{abstract}
We prove a Lichnerowicz type lower bound for the first nontrivial eigenvalue of the $p$-Laplacian on K\"ahler manifolds.  Parallel to the $p=2$ case, the first eigenvalue lower bound is improved by using a decomposition of the Hessian on K\"ahler manifolds with positive Ricci curvature.
\end{abstract}
\section{Introduction}
Let $(M,g)$ be a $n$-dimensional compact Riemannian manifold, possibly with boundary.  The $p$-Laplace operator $\Delta_p$ is defined by
\begin{equation*}
\Delta_p(f) := \diver(|\nabla f|^{p-2}\nabla f).
\end{equation*}
This is a generalization of the classical Laplace operator ($p=2$) and has found many applications in mathematics as well as physics.   While it is only a quasilinear elliptic operator for $p\neq 2$, the $p$-Laplacian shares many characteristics to the classical Laplacian.  See, for instance, \cite{lindqvist2}, \cite{lindqvist} for a general reference on the $p$-Laplacian.
The corresponding $p$-Laplace eigenvalue equation is given by
\begin{equation*}
\Delta_p(f) = -\mu|f|^{p-2}f,
\end{equation*}
with appropriate boundary conditions.  This equation arises from the following variational characterization of the first nonzero eigenvalue given by
\begin{equation*}
\mu_{1,p} = \inf\left\{ \frac{\int_M|\nabla f|^p}{\int_M|f|^p} \ | \ f \in W^{1,p}(M)\backslash\{0\}, \int_M |f|^{p-2}f=0\right\}
\end{equation*}
for closed $M$ and
\begin{equation*}
\lambda_{1,p} = \inf\left\{ \frac{\int_M|\nabla f|^p}{\int_M|f|^p} \ | \ f \in W_c^{1,p}(M)\backslash\{0\}\right\}
\end{equation*}
if we impose the Dirichlet boundary condition.  Note that unlike the case $p=2$, the eigenfunctions have only partial regularity, i.e., of class $C^{1,\alpha}$ and for $\mu_{1,p} \neq 0$, they are never $C^2$ (c.f. \cite{kawai-nakauchi}).  Note that $f$ is smooth away from the set $\{\nabla f=0\}$.
In \cite{seto-wei}, a Lichnerowicz-type lower bound was established for $\mu_{1,p}$, namely, on complete $n$-dimensional Riemannian manifolds with $\Ric \geq Kg$, $K>0$, and $p \geq 2$,
\begin{equation*}
\mu_{1,p}^{\frac{2}{p}} \geq \left(1+ \frac{1}{\sqrt{n}(p-2)+n-1}\right)\frac{K}{p-1}.
\end{equation*}
In fact, this was shown in a slightly more general context of integral Ricci curvature conditions.   Here we show that the lower bound can be improved on K\"ahler manifolds.
\begin{theorem}\label{mainthm}
Let $(M,J,g)$ be an $n=2m$ (real) dimensional K\"ahler manifold, possibly with boundary.  Assume that the underlying (real) Ricci curvature satisfies $\Ric \geq Kg$ for some constant $K>0$. If $\dd M = \emptyset$, then for $p\geq 2$,
\begin{equation}\label{dirichlet}
\mu_{1,p}^{\frac{2}{p}} \geq \frac{p+2}{p(p-1)}K = \left( 1+ \frac{2}{p}\right)\frac{K}{p-1}.
\end{equation}
If $\dd M \neq \emptyset$, we assume the convexity condition that $\displaystyle\frac{p}{2}H+\SecFun(J\mathbf{n},J\mathbf{n})\geq 0$ and the Dirichlet boundary condition, where $\mathbf{n}$ is the unit outward normal vector field on $\dd M$, $H$ is the mean curvature, and $\SecFun$ is the second fundamental form.  Then for $p\geq 2$, 
\begin{equation}\label{neumann}
\lambda_{1,p}^{\frac{2}{p}} \geq \frac{p+2}{p(p-1)}K.
\end{equation}
\end{theorem}
When $p=2$, this recovers the results of Urakawa \cite{urakawa} for the closed case and Guedj, Kolev, and Yeganefar \cite{guedj-kolev-yeganefar} for the Dirichlet boundary case.  See also \cite{grosjean} and \cite{lichnerowicz} regarding the lower bound when $p=2$.  For upper bounds, Chen and Wei \cite{chen-wei} provide some estimates for the $p$-Laplacian on submanifolds of space forms.

To obtain our estimate, we first establish a Reilly type formula for the $p$-Laplacian. The main difficulty for the $p>2$ case is the introduction of the term involving an inner product of the Hessian in the $\nabla f$ direction with the same term but pushed forward by the complex structure $J$.  As there is no a priori relation between the eigenfunction $f$ with the complex structure $J$, unlike the Riemannian case, we need to take advantage of all terms involved in the $p$-Bochner formula.  

\begin{remark}
Using the methods of \cite{seto-wei}, we can show for $p>2$ that a lower bound holds under the assumption of integral Ricci curvature.  See Remark \ref{integral}.
\end{remark}

In \S 2, we give some backgrounds concerning manifolds with boundary and give a Reilly formula adapted for the $p$-Laplacian case.  In \S 3, we give some detail for the decomposition of the Hessian on K\"ahler manifolds and prove the eigenvalue lower bound by applying this decomposition to the Reilly formula.  

\begin{acknowledgements}
The authors would like to thank Professor Guofang Wei for her interest and valuable comments on the initial draft, as well as a reference to the Reilly-type formulas.  The authors would also like to thank the referee whose careful proofreading and comments have greatly improved the paper.
\end{acknowledgements}

\section{$p$-Reilly formula}
Let $(M,g)$ be a compact Riemannian manifold with boundary.
\begin{definition}
The second fundamental form is
\begin{align*}
\SecFun(X,Y) = \langle \nabla_X\mathbf{n},Y\rangle,
\end{align*}
where $\mathbf{n}$ is the unit outward normal vector on $\dd M$.
\end{definition}
We begin with the following basic fact.
\begin{lemma}[(8.1) \cite{li}]
Let $S^m\subset N^n$ be an $m$-dimensional submanifold of an arbitrary manifold $N$ and let $\{e_i\}_{i=1}^m$ be an adapted orthonormal frame tangential to $S$ and $\{e_{\nu}\}_{\nu=m+1}^{n}$ normal to $S$.  Then for $1\leq i,j \leq m$, the Hessian is related by
\begin{equation*}
(\Hess_N f)_{ij} =(\Hess_S f)_{ij}+\sum_{\nu=m+1}^n \SecFun_{ij} e_{\nu}f.
\end{equation*}
\end{lemma}
Specializing to hypersurfaces $\bar M^{n-1}\subset M^{n}$, we take the trace to get
\begin{equation}\label{laplacehypersurf}
\Delta f - f_{nn}= \Delta_{\bar{M}} f + H\frac{\dd f}{\dd n},
\end{equation}
where $H$ is the mean curvature and $\Delta_{\bar{M}}$ is the Laplacian on $\bar{M}^{n-1}$.   

As noted in \cite{guedj-kolev-yeganefar}, on K\"ahler manifolds, we have the following decomposition of the Hessian into the sum of a $J$-symmetric bilinear form and a $J$-skew-symmetric bilinear form:
\begin{align*}
\Hess f = H_1f +H_2f
\end{align*}
where
\begin{align*}
H_1f(X,Y) &= \frac{1}{2}(\Hess f(X,Y)+\Hess f(JX,JY))\\
H_2f(X,Y) &= \frac{1}{2}(\Hess f(X,Y) - \Hess f(JX,JY)).
\end{align*}
Here the skew-symmetrization of $H_1$ will lead to the $(1,1)$-Hessian and $H_2$ is the $(2,0)+(0,2)$ Hessian.  Under this decomposition,
\begin{align*}
2\|H_1 f\|^2 &= \|\Hess f\|^2 + \langle \Hess f, J^*\Hess f\rangle\\
2\|H_2 f\|^2 &= \|\Hess f\|^2 -\langle \Hess f,J^*\Hess f\rangle.
\end{align*}
Note that the above holds for complex manifolds and does not require that the complex structure be covariantly constant.  The K\"ahler structure is used later when we want to relate $\langle \Hess f,J^*\Hess f\rangle$ to a curvature term.

We first establish a $p$-Reilly formula, 
\begin{lemma}[$p$-Reilly formula]
For $f\in C^2(M)$ and $p\geq 2$, 
\begin{align}
\begin{split}\label{reillyform}
&\int_{\dd M}|\nabla f|^{p-2}\left\{ -(\Delta_{\dd M} f+H\nabla_nf)\nabla_nf-\SecFun(\nabla_{\dd M}f,\nabla_{\dd M}f) + \langle \nabla (\nabla_n f),\nabla f\rangle_{\dd M} \right\} \\
&=(p-2)\int_M|\nabla f|^{p-2}|\nabla|\nabla f||^2  - \int_M(\Delta f)(\Delta_p f) \\
&\hspace{0.2 in} + \int_M|\nabla f|^{p-2}(2|H_2f|^2+\Ric(\nabla f,\nabla f)+\langle \Hess f,J^*\Hess f\rangle ).
\end{split}
\end{align}
\end{lemma}
\begin{remark}
See also a related Reilly type formula on K\"ahler manifolds in \cite{wang}, and a similar $p$-Reilly formula in \cite{wang-li}.  Here we used the decomposition of the Hessian using $H_2$.  If instead we use the decomposition with $H_1$, then we would obtain a Reilly formula similar to the one presented in \cite{wang}, where for $p=2$, the Ricci term cancels out.  Since we want to take advantage of the Ricci curvature lower bound, this version is not suitable for our application.
\end{remark}
\begin{proof}
We integrate the following $p$-Bochner formula (Lemma 3.1 \cite{seto-wei}, note the typo in the statement there but is otherwise used correctly in its application).
\begin{equation*}
\frac{1}{p}\Delta(|\nabla f|^p) = (p-2)|\nabla f|^{p-2}|\nabla|\nabla f||^2+|\nabla f|^{p-2}\left\{|\Hess f|^2 +\langle \nabla f,\nabla \Delta f\rangle + \Ric(\nabla f,\nabla f)\right\}.
\end{equation*}
Integrating the left hand side, we have
\begin{align*}
\frac{1}{p}\int_M \Delta(|\nabla f|^p) &= \frac{1}{p}\int_{\dd M} \nabla_n|\nabla f|^pdS \\
&= \int_{\dd M}|\nabla f|^{p-2}\langle \nabla_n\nabla f,\nabla f\rangle.
\end{align*}
Pointwise, using an (adapted) orthonormal frame $\{e_i\}$ with $e_n = \mathbf{n}$ and \eqref{laplacehypersurf} we have
\begin{align*}
\langle \nabla_n \nabla f, \nabla f\rangle 
&=\Hess f(e_n,e_n)\nabla_n f+\sum_{i=1}^{n-1} \Hess f(e_n,e_i)\nabla_i f\\
&=(\Delta f - \Delta_{\dd M} f -H\nabla_n f)\nabla_n f + \sum_{i=1}^{n-1} \Hess f(e_n,e_i)\nabla_i f.
\end{align*}
For fixed $i \leq n-1$, we have
\begin{align*}
\Hess f(e_n,e_i)
&=\sum_{j=1}^{n-1}\langle \nabla_i (\nabla_j fe_j),e_n\rangle + \langle \nabla_i (\nabla_nf e_n),e_n\rangle\\
&=-\sum_{j=1}^{n-1}\langle \nabla_j fe_j,\nabla_i e_n\rangle + e_i(\nabla_nf)-\nabla_nf \langle e_n,\nabla_ie_n\rangle\\
&=-\sum_{j=1}^{n-1}(\nabla_jf)\langle \nabla_ie_n,e_j\rangle + e_i(\nabla_n f)\\
&=-\sum_{j=1}^{n-1}\SecFun_{ij}(\nabla_jf) + e_i(\nabla_n f).
\end{align*}
Combining the above equations, we get
\begin{align}
\begin{split}\label{boundaryintegral}
&\int_{\dd M} |\nabla f|^{p-2}\langle \nabla_n\nabla f,\nabla f\rangle\\ 
&=\int_{\dd M}|\nabla f|^{p-2}\left\{(\Delta f)\nabla_n f-(\Delta_{\dd M} f) \nabla_nf-H(\nabla_n f)^2-\SecFun(\nabla_{\dd M} f,\nabla_{\dd M} f) + \langle \nabla(\nabla_n f),\nabla f\rangle_{\dd M}\right\}.
\end{split}
\end{align}
Integrating the right hand side of the $p$-Bochner formula, for the third term we integrate by parts to obtain
\begin{align*}
\int_M|\nabla f|^{p-2}\langle \nabla f, \nabla \Delta f\rangle &= \int_M \diver(|\nabla f|^{p-2}(\Delta f)\nabla f)-\int_M\Delta f\Delta_p f\\
&=\int_{\dd M} \nabla_nf|\nabla f|^{p-2}\Delta f - \int_M\Delta f\Delta_pf.
\end{align*}
Using the decomposition of the Hessian,
\begin{align*}
\int_M|\nabla f|^{p-2}|\Hess f|^2 &= \int_M2|\nabla f|^{p-2}|H_2f|^2 + |\nabla f|^{p-2}\langle \Hess f, J^*\Hess f\rangle
\end{align*}
and combining the equations, we obtain the result. 
\end{proof}

\section{Proof of Theorem \ref{mainthm}}
To obtain the Lichnerowicz estimate for $p=2$, one usually applies the Cauchy-Schwarz inequality to the norm of the Hessian to relate to the Laplacian.  On K\"ahler manifolds, we can take advantage of the decomposition of the Hessian which contains a curvature term.  This was a key observation in \cite{guedj-kolev-yeganefar} and we modify to the $p$-Laplacian case.  Consider the term
\begin{align}
\begin{split}\label{mainprodrule}
\diver(|\nabla f|^{p-2}J^*\Hess f(\nabla f, \cdot)^{\#}) &= \langle \nabla |\nabla f|^{p-2},J^*\Hess f(\nabla f, \cdot)^{\#}\rangle + |\nabla f|^{p-2}\diver(J^*\Hess f(\nabla f, \cdot )^{\#}).
\end{split}
\end{align}
Using an (adapted) orthonormal frame $\{e_i\}$ with $e_n = \mathbf{n}$, the second term on the right hand side of \eqref{mainprodrule} is expressed locally as
\begin{align}
\begin{split}\label{prodrule}
\diver(\Hess f(J\nabla f, J\cdot)^{\#}) &= \sum_{i=1}^n e_i\langle \nabla_{Je_i}\nabla f,J\nabla f\rangle \\
&=\sum_{i= 1}^n\langle \nabla_{e_i}\nabla_{Je_i}\nabla f,J\nabla f\rangle +\langle \nabla_{Je_i}\nabla f,J\nabla_{e_i}\nabla f\rangle.
\end{split}
\end{align}
Here we used the fact that $\nabla J =0$.  The first term on the right hand side of \eqref{prodrule} can be modified in the following way:  We are tracing over an orthonormal frame $\{e_i\}$, so instead, we trace over the frame $\{Je_i\}$.  Then
\begin{align*}
\sum_{i=1}^n \langle \nabla_{e_i}\nabla_{Je_i}\nabla f,J\nabla f\rangle &= \frac{1}{2}\sum_{i=1}^n \langle \nabla_{e_i}\nabla_{Je_i}\nabla f,J\nabla f\rangle - \langle \nabla_{Je_i}\nabla_{e_i}\nabla f,J\nabla f\rangle\\
&=\frac{1}{2}\sum_{i=1}^n\langle (\nabla_{e_i}\nabla_{Je_i}-\nabla_{Je_i}\nabla_{e_i})\nabla f, J\nabla f\rangle \\
&=-\frac{1}{2}\sum_{i=1}^nR(e_i,Je_i,\nabla f,J\nabla f\rangle \\
&=-\frac{1}{2}\sum_{i=1}^n R(e_i,\nabla f, e_i,\nabla f)+R(e_i, J\nabla f, e_i, J\nabla f) \\
&=-\Ric(\nabla f, \nabla f),
\end{align*}
where the second to last line uses the Bianchi identity.  The second term on the right hand side of \eqref{prodrule} is given locally as
\begin{align*}
\sum_{i=1}^n \langle \nabla_{Je_i}\nabla f,J\nabla_{e_i}\nabla f\rangle &= -\sum_{i=1}^n \langle J\nabla_{Je_i}\nabla f, \nabla_{e_i}\nabla f\rangle \\
&=-\sum_{i,j=1}^n\langle \langle J\nabla_{Je_i}\nabla f,e_j\rangle e_j,\nabla_{e_i}\nabla f\rangle \\
&= \sum_{i,j=1}^n\langle \nabla_{e_i}\nabla f,e_j\rangle\langle \nabla_{Je_i}\nabla f,Je_j\rangle \\
&= \langle \Hess f, J^*\Hess f\rangle.
\end{align*}
For the first term on the right hand side of \eqref{mainprodrule} we can rewrite as
\begin{align*}
\langle \nabla |\nabla f|^{p-2},\Hess f(J\nabla f, J\cdot)^\#\rangle &= (p-2)|\nabla f|^{p-4}\langle \nabla_{Je_i}\nabla f, J\nabla f\rangle \Hess f(\nabla f, e_i)\\
&=(p-2)|\nabla f|^{p-4}\langle \nabla_{Je_i}\nabla f,J\nabla f\rangle\langle\nabla_{e_i}\nabla f,\nabla f\rangle\\
&=-(p-2)|\nabla f|^{p-4}\langle \nabla f, e_j\rangle\langle\nabla f,e_k\rangle\langle J\nabla_{Je_i}\nabla f, e_j\rangle \langle \nabla_{e_i}\nabla f,e_k\rangle\\
&=(p-2)|\nabla f|^{p-4}\langle \Hess f(\nabla f,\cdot),J^*\Hess f(\nabla f,\cdot)\rangle.
\end{align*}
Combining the above equations, we get
\begin{align*}
\diver(|\nabla f|^{p-2}J^*\Hess f(\nabla f, \cdot)^{\#}) &= -|\nabla f|^{p-2}\Ric(\nabla f,\nabla f) + |\nabla f|^{p-2}\langle \Hess f,J^*\Hess f\rangle\\
&\hspace{0.2 in}+(p-2)|\nabla f|^{p-4}\langle \Hess f(\nabla f,\cdot),J^*\Hess f(\nabla f, \cdot)\rangle.
\end{align*}
Applying divergence theorem to the above equation, the integrand of the boundary term is
\begin{align*}
|\nabla f|^{p-2}J^*\Hess f(\nabla f,  e_n) &= |\nabla f|^{p-2}J^*\Hess f(\nabla_{\dd M}f, e_n) + |\nabla f|^{p-2}(\nabla_n f)J^*\Hess f(e_n,e_n).
\end{align*}
From the decomposition
\begin{align*}
\nabla_XY &= \sum_{i=1}^{n-1}\langle \nabla_XY,e_i\rangle e_i + \langle \nabla_XY,n\rangle n \\
&=(\nabla_X)_{\dd M}Y - \SecFun(X,Y)n,
\end{align*}
for $X,Y \in T_p(\dd M)$ and
\begin{align*}
\Hess f(X,Y) = \Hess f_{\dd M}(X,Y) + (\nabla_n f)\SecFun(X,Y)
\end{align*}
we have
\begin{align*}
|\nabla f|^{p-2}J^*\Hess f(\nabla f,  e_n) &= |\nabla f|^{p-2}J^*\Hess f(\nabla_{\dd M}f, e_n) + |\nabla f|^{p-2}(\nabla_n f)\Hess f_{\dd M}(Je_n,Je_n) \\
&\hspace{0.2 in} + |\nabla f|^{p-2}(\nabla_n f)^2\SecFun(Je_n,Je_n).
\end{align*}
Therefore,
\begin{align}
\begin{split}\label{hessterm}
\int_M&|\nabla f|^{p-2}\langle \Hess f, J^*\Hess f\rangle +(p-2)\int_M |\nabla f|^{p-4} \langle \Hess f(\nabla f,\cdot),J^*\Hess f(\nabla f,\cdot)\rangle\\
 &= \int_M|\nabla f|^{p-2}\Ric(\nabla f,\nabla f)+\int_{\dd M}|\nabla f|^{p-2}J^*\Hess f(\nabla_{\dd M}f,e_n)\\
&\hspace{0.2 in}+\int_{\dd M}|\nabla f|^{p-2}(\nabla_n f)\Hess f_{\dd M}(Je_n,Je_n)+\int_{\dd M}|\nabla f|^{p-2}(\nabla_n f)^2\SecFun(Je_n,Je_n).
\end{split}
\end{align}
Combining \eqref{hessterm} with the Reilly formula \eqref{reillyform},
\begin{align}
\begin{split}\label{keyineq}
&\int_{\dd M}|\nabla f|^{p-2}\left\{ -(\Delta_{\dd M} f+H\nabla_nf)\nabla_nf-\SecFun(\nabla_{\dd M}f,\nabla_{\dd M}f) + \langle \nabla (\nabla_n f),\nabla f\rangle_{\dd M} \right\} \\
&=(p-2)\int_M|\nabla f|^{p-2}|\nabla|\nabla f||^2  - \int_M(\Delta f)(\Delta_p f) \\
&\hspace{0.2 in} + \int_M|\nabla f|^{p-2}(2|H_2f|^2+2\Ric(\nabla f,\nabla f))\\
&\hspace{0.2 in}-(p-2)\int_M|\nabla f|^{p-4}\langle \Hess f(\nabla f,\cdot),J^*\Hess f(\nabla f,\cdot)\rangle \\
&\hspace{0.2 in} +\int_{\dd M}|\nabla f|^{p-2}J^*\Hess f(\nabla_{\dd M}f,e_n)+\int_{\dd M}|\nabla f|^{p-2}(\nabla_n f)\Hess f_{\dd M}(Je_n,Je_n)\\
&\hspace{0.2 in}+\int_{\dd M}|\nabla f|^{p-2}(\nabla_n f)^2\SecFun(Je_n,Je_n).
\end{split}
\end{align}
Since
\begin{align*}
|\nabla |\nabla f||^2 = |\Hess f(\nabla f,\cdot)|^2|\nabla f|^{-2},
\end{align*}
we can use the decomposition of the Hessian so that
\begin{align*}
\int_M |\nabla f|^{p-2}|\nabla |\nabla f||^2 &= \int_M |\nabla f|^{p-4}|\Hess f(\nabla f,\cdot)|^2 \\
&= \int_M |\nabla f|^{p-4}( 4|H_2 f(\nabla f,\cdot)|^2 - |\Hess f(J\nabla f, J\cdot)|^2)\\
&\hspace{0.2 in}+2\int_M |\nabla f|^{p-4}\langle \Hess f(\nabla f,\cdot),\Hess f(J\nabla f,J\cdot)\rangle) \\
&\geq \int_M |\nabla f|^{p-4}( 4|H_2 f(\nabla f,\cdot)|^2 -\int_M |\nabla f|^{p-2}|\Hess f|^2\\
&\hspace{0.2 in}+2\int_M |\nabla f|^{p-4}\langle \Hess f(\nabla f,\cdot),\Hess f(J\nabla f,J\cdot)\rangle) .
\end{align*}
The $|\Hess f|^2$ term can be rewritten as
\begin{align*}
-&\int_M|\nabla f|^{p-2}|\Hess f|^2\\
 &= -\int_M|\nabla f|^{p-2}\diver(\Hess f(\nabla f,\cdot)) + \int_M |\nabla f|^{p-2}\langle \Delta \nabla f,\nabla f\rangle\\
&=-\int_M\diver(|\nabla f|^{p-2}\Hess f(\nabla f,\cdot)) + \int_Me_i(|\nabla f|^{p-2})\Hess f(\nabla f, e_i) + \int_M |\nabla f|^{p-2} \langle \Delta \nabla f,\nabla f\rangle\\
&=-\int_M\diver(|\nabla f|^{p-2}\Hess f(\nabla f,\cdot)) +(p-2)\int_M|\nabla f|^{p-4}|\Hess f(\nabla f,\cdot)|^2+ \int_M |\nabla f|^{p-2} \langle \Delta \nabla f,\nabla f\rangle.
\end{align*}
The last term can be written in terms of the $p$-Laplacian as
\begin{align*}
\int_M |\nabla f|^{p-2}\langle \Delta \nabla f,\nabla f\rangle &= \int_M|\nabla f|^{p-2}\Ric(\nabla f,\nabla f) + \int_M|\nabla f|^{p-2}\langle \nabla (\Delta f),\nabla f\rangle \\
&= \int_M|\nabla f|^{p-2}\Ric(\nabla f,\nabla f)-\int_M \Delta f \Delta_p f + \int_{\dd M}\nabla_n f|\nabla f|^{p-2}\Delta f.
\end{align*}
Combining these together and dropping the non-negative terms, we have for $p \geq 2$, 
\begin{align*}
\frac{(p-2)}{2}\int_M & |\nabla f|^{p-2}|\nabla |\nabla f||^2 \\&\geq (p-2)\int_M |\nabla f|^{p-4}\langle \Hess f(\nabla f, \cdot),\Hess f(J\nabla f, J\cdot)\rangle\\
&\hspace{0.2 in} + \frac{(p-2)}{2}\int_M|\nabla f|^{p-2}\Ric(\nabla f,\nabla f) - \frac{(p-2)}{2}\int_M\Delta f \Delta_p f\\
&\hspace{0.2 in}+ \frac{(p-2)}{2}\int_{\dd M}\nabla_n f|\nabla f|^{p-2}\Delta f -\frac{(p-2)}{2} \int_{\dd M}|\nabla f|^{p-2}\Hess f(\nabla f, n).
\end{align*}
The boundary term can be simplified using \eqref{boundaryintegral} so that
\begin{align*}
&\frac{(p-2)}{2}\int_{\dd M}|\nabla f|^{p-2}((\Delta f)\nabla_n f - \langle \nabla_n\nabla f,\nabla f\rangle) \\
&=\frac{(p-2)}{2}\int_{\dd M}|\nabla f|^{p-2}\left\{((\Delta_{\dd M} f)+H\nabla_n f)\nabla_n f+\SecFun(\nabla_{\dd M}f,\nabla_{\dd M}f)-\langle \nabla(\nabla_n f),\nabla f\rangle_{\dd M}\right\}.
\end{align*}
Combining the above with \eqref{keyineq}, we get
\begin{align}
\begin{split}\label{mainform}
&\frac{p}{2}\int_{\dd M}|\nabla f|^{p-2}\left\{ -(\Delta_{\dd M} f+H\nabla_nf)\nabla_nf-\SecFun(\nabla_{\dd M}f,\nabla_{\dd M}f) + \langle \nabla (\nabla_n f),\nabla f\rangle_{\dd M} \right\} \\
&\geq - \frac{p}{2}\int_M(\Delta f)(\Delta_p f)+ \frac{(p+2)}{2}\int_M|\nabla f|^{p-2}\Ric(\nabla f,\nabla f) \\
&\hspace{0.2 in} +\int_{\dd M}|\nabla f|^{p-2}J^*\Hess f(\nabla_{\dd M}f,e_n)+\int_{\dd M}|\nabla f|^{p-2}(\nabla_n f)\Hess f_{\dd M}(Je_n,Je_n)\\
&\hspace{0.2 in}+\int_{\dd M}|\nabla f|^{p-2}(\nabla_n f)^2\SecFun(Je_n,Je_n).
\end{split}
\end{align}
Now we are ready to prove Theorem \ref{mainthm}.
\begin{proof}
By a density argument, we can apply \eqref{mainform} to the first eigenfunction $f$ and in particular, for $\Ric \geq K$, 
\begin{align*}
\frac{(p+2)}{2}\int_M|\nabla f|^{p-2}\Ric(\nabla f,\nabla f) &\geq \frac{(p+2)K}{2}\int_M |\nabla f|^p = \frac{(p+2)K}{2}\lambda_{1,p}\int_M|f|^p
\end{align*}
and
\begin{align*}
-\frac{p}{2}\int_M (\Delta f)(\Delta_p f) &= \frac{p}{2}\lambda_{1,p}\int_M|f|^{p-2}f\Delta f\\
&=-\frac{p}{2}\lambda_{1,p}\int_M\langle \nabla(|f|^{p-2}f),\nabla f\rangle \\
&=-\frac{p(p-1)}{2}\lambda_{1,p}\int_M |f|^{p-2}|\nabla f|^2 \\
&\geq -\frac{p(p-1)}{2}\lambda_{1,p}\left(\int_M |f|^p\right)^{1-\frac{2}{p}}\left(\int_M |\nabla f|^p\right)^{\frac{2}{p}}\\
&=-\frac{p(p-1)}{2}\lambda_{1,p}^{1+\frac{2}{p}}\int_M |f|^p.
\end{align*}

Using Dirichlet boundary condition and the above inequalities \eqref{mainform} becomes
\begin{align*}
-&\frac{p}{2}\int_{\dd M}H|\nabla f|^{p-2}(\nabla_n f)^2 \\
&\geq \left(\frac{(p+2)K}{2}\lambda_{1,p}-\frac{p(p-1)}{2}\lambda_{1,p}^{1+\frac{2}{p}}\right)\int_M |f|^p + \int_{\dd M}|\nabla f|^{p-2}(\nabla_n f)^2\SecFun(Je_n,Je_n).
\end{align*}
Therefore,
\begin{align*}
\frac{\lambda_{1,p}}{2}\left(\lambda_{1,p}^{\frac{2}{p}}p(p-1)-(p+2)K\right)\int_M|f|^p \geq \int_{\dd M}\left(\frac{p}{2}H+\SecFun(Je_n,Je_n)\right)|\nabla f|^{p-2}(\nabla_n f)^2.
\end{align*}
By the convexity condition, the expression must be nonnegative therefore
\begin{equation*}
\lambda_{1,p}^{\frac{2}{p}}\geq \frac{p+2}{p(p-1)}K.
\end{equation*}
The same conclusion holds for $\mu_{1,p}$ since the boundary integrals are zero in this case.
\end{proof}

\begin{remark}\label{integral}
By following the methods used in \cite{seto-wei}, when $p>2$, one can use the remaining term $\frac{(p-2)}{2}|\nabla|\nabla f||^2$ which we dropped to obtain a lower bound under integral Ricci curvature condition as well.  In detail, for each $x\in M$, let  $\rho\left( x\right) $ denote the smallest
eigenvalue for the Ricci tensor $\mathrm{Ric}:T_{x}M\rightarrow
T_{x}M,$ and $\Ric_-^K(x) = \left( (n-1)K - \rho
(x)\right)_+ = \max \left\{ 0, (n-1)K - \rho (x) \right\}$, the amount of Ricci curvature lying below $(n-1)K$. Let
\begin{equation*} \| \Ric_-^K \|^*_{q, R} = \sup_{x\in M} \left( \frac{1}{\vol(B(x,R))}\int_{B\left( x,R\right) } (\mathrm{Ric}_-^K)^{q}\,
d vol\right)^{\frac 1q}.
\end{equation*}
Then  $\| \mathrm{Ric}_-^K \|^*_{q,R} $ measures the amount of Ricci
curvature lying below a given bound, in this case,  $(n-1)K$, in the $L^q$ sense.  Then for a complete manifold $M$ with $q>\frac{n}{2}$, $p\geq 2$ and $K>0$, there exists $\eps = \eps(n,p,q,K)$ such that if $\|\Ric_-^K\|_q^* \leq \eps$, then
\begin{align*}
\mu_{1,p}^{\frac{2}{p}} \geq \left(1+\frac{2}{p}\right)\left(\frac{K}{p-1}-\frac{2}{p-1}\|\Ric_-^K\|_q^*\right).
\end{align*}
\end{remark}


\begin{bibdiv}
\begin{biblist}
\bib{chen-wei}{article}{
   author={Chen, Hang},
   author={Wei, Guofang},
   title={Reilly-type inequalities for $p$-Laplacian on submanifolds in space forms},
   journal={arXiv:1806.09061},
   date={2018}
}

\bib{grosjean}{article}{
   author={Grosjean, Jean-Fran\c{c}ois},
   title={A new Lichnerowicz-Obata estimate in the presence of a parallel
   $p$-form},
   journal={Manuscripta Math.},
   volume={107},
   date={2002},
   number={4},
   pages={503--520},
   issn={0025-2611},
   review={\MR{1906773}},
   doi={10.1007/s002290200248},
}

\bib{guedj-kolev-yeganefar}{article}{
   author={Guedj, Vincent},
   author={Kolev, Boris},
   author={Yeganefar, Nader},
   title={A Lichnerowicz estimate for the first eigenvalue of convex domains
   in K\"ahler manifolds},
   journal={Anal. PDE},
   volume={6},
   date={2013},
   number={5},
   pages={1001--1012},
   issn={2157-5045},
   review={\MR{3125547}},
   doi={10.2140/apde.2013.6.1001},
}

\bib{kawai-nakauchi}{article}{
   author={Kawai, Shigeo},
   author={Nakauchi, Nobumitsu},
   title={The first eigenvalue of the $p$-Laplacian on a compact Riemannian
   manifold},
   journal={Nonlinear Anal.},
   volume={55},
   date={2003},
   number={1-2},
   pages={33--46},
   issn={0362-546X},
   review={\MR{2001630}},
   doi={10.1016/S0362-546X(03)00209-8},
}

\bib{li}{book}{
   author={Li, Peter},
   title={Geometric analysis},
   series={Cambridge Studies in Advanced Mathematics},
   volume={134},
   publisher={Cambridge University Press, Cambridge},
   date={2012},
   pages={x+406},
   isbn={978-1-107-02064-1},
   review={\MR{2962229}},
   doi={10.1017/CBO9781139105798},
}

\bib{lichnerowicz}{book}{
   author={Lichnerowicz, Andr\'e},
   title={G\'eom\'etrie des groupes de transformations},
   language={French},
   publisher={Travaux et Recherches Math\'ematiques, III. Dunod, Paris},
   date={1958},
   pages={ix+193},
   review={\MR{0124009}},
}

\bib{lindqvist2}{article}{
   author={Lindqvist, Peter},
   title={A nonlinear eigenvalue problem},
   conference={
      title={Topics in mathematical analysis},
   },
   book={
      series={Ser. Anal. Appl. Comput.},
      volume={3},
      publisher={World Sci. Publ., Hackensack, NJ},
   },
   date={2008},
   pages={175--203},
   review={\MR{2462954}},
   doi={10.1142/97898128110660005},
}

\bib{lindqvist}{book}{
   author={Lindqvist, Peter},
   title={Notes on the $p$-Laplace equation},
   series={Report. University of Jyv\"askyl\"a Department of Mathematics and
   Statistics},
   volume={102},
   publisher={University of Jyv\"askyl\"a, Jyv\"askyl\"a},
   date={2006},
   pages={ii+80},
   isbn={951-39-2586-2},
   review={\MR{2242021}},
}

\bib{reilly}{article}{
   author={Reilly, Robert C.},
   title={Applications of the Hessian operator in a Riemannian manifold},
   journal={Indiana Univ. Math. J.},
   volume={26},
   date={1977},
   number={3},
   pages={459--472},
   issn={0022-2518},
   review={\MR{0474149}},
   doi={10.1512/iumj.1977.26.26036},
}

\bib{seto-wei}{article}{
   author={Seto, Shoo},
   author={Wei, Guofang},
   title={First eigenvalue of the $p$-Laplacian under integral curvature
   condition},
   journal={Nonlinear Anal.},
   volume={163},
   date={2017},
   pages={60--70},
   issn={0362-546X},
   review={\MR{3695968}},
   doi={10.1016/j.na.2017.07.007},
}

\bib{urakawa}{article}{
   author={Urakawa, Hajime},
   title={Stability of harmonic maps and eigenvalues of the Laplacian},
   journal={Trans. Amer. Math. Soc.},
   volume={301},
   date={1987},
   number={2},
   pages={557--589},
   issn={0002-9947},
   review={\MR{882704}},
   doi={10.2307/2000659},
}

\bib{wang}{article}{
   author={Wang, Xiaodong},
   title={An integral formula in K\"ahler geometry with applications},
   journal={Commun. Contemp. Math.},
   volume={19},
   date={2017},
   number={5},
   pages={1650063, 12},
   issn={0219-1997},
   review={\MR{3670795}},
   doi={10.1142/S0219199716500632},
}

\bib{wang-li}{article}{
   author={Wang, Yu-Zhao},
   author={Li, Huai-Qian},
   title={Lower bound estimates for the first eigenvalue of the weighted
   $p$-Laplacian on smooth metric measure spaces},
   journal={Differential Geom. Appl.},
   volume={45},
   date={2016},
   pages={23--42},
   issn={0926-2245},
   review={\MR{3457386}},
   doi={10.1016/j.difgeo.2015.11.008},
}

\end{biblist}
\end{bibdiv}
\end{document}